\documentclass[a4paper,10pt]{article}

\usepackage[english]{babel}
\usepackage[utf8x]{inputenc}
\usepackage{amsmath}
\usepackage{tikz}
\usepackage{color}
\usepackage{graphicx}
\usepackage[colorinlistoftodos]{todonotes}
\usepackage{amssymb}
\usepackage{amsfonts}
\usepackage{amsthm}
\usepackage{mathrsfs}
\usepackage{url}
\usepackage{setspace}
\newtheorem{thm}{Theorem}[section]

\newtheorem{cnj}[thm]{Conjecture}

\newtheorem{cor}[thm]{Corollary}

\def\a{{\alpha}}

\def\cF{{\cal F}}
\def\cI{{\cal I}}
\def\cL{{\cal L}}

\def\pr{{\prime}}
\def\sse{{\subseteq}}

\def\flip{{\sf flip}}
\def\shift{{\sf shift}}
\def\slide{{\sf slide}}

\definecolor{brwn}{RGB}{140, 70, 20}
\definecolor{gren}{RGB}{  0,140, 10}

\author{
Glenn Hurlbert\thanks{
Department of Mathematics and Applied Mathematics,
Virginia Commonwealth University, 
Richmond, VA, USA, 
\texttt{ghurlbert@vcu.edu}.
}\\
Vikram Kamat\thanks{
Department of Mathematics \& Statistics, 
Villanova University, 
Villanova, PA, USA, 
\texttt{vikram.kamat@villanova.edu}}
}

\title{On intersecting families of independent sets in trees}
\date{}

\begin{document}
\maketitle

\begin{abstract}
A family of sets is \textit{intersecting} if every pair of its sets intersect.
A star is a family with some element (a {\it center}) in each of its sets.
The classical result of Erd\H{o}s, Ko, and Rado (1961) states that every intersecting family of $r$-subsets of $[n]$ with $r\leq n/2$ has size at most that of a star. 
Let $G$ be a graph, $\a(G)$ be its independence number, and $\mu(G)$ be the size of a smallest maximal independent set in $G$. 
We say that $G$ is $r$-EKR if, among all maximum-sized intersecting families of independent $r$-subsets of vertices of $G$, there is a star.
In 2005 Holroyd and Talbot conjectured that every graph $G$ is $r$-EKR for all $1\leq r\leq \mu(G)/2$. 
We verified the conjecture in 2011 for all chordal graphs containing an isolated vertex.

For a graph without isolated vertices it is difficult to determine a center of a largest star, which is often necessary to prove that it is EKR.
A tree has the {\it leaf property} if some largest star occurs on one of its leaves.
We proved that every tree $T$ has the leaf property when $r\le 4$, and in 2017 Borg and other authors gave examples of families of trees not having the leaf property when $r\ge 5$.
A {\it split} vertex in a tree is a vertex of degree at least 3.
A {\it spider} is a tree with exactly one split vertex.
Here we prove that all spiders $T$ have the leaf property for all $r\le \a(T)$, and we characterize which of its leaves is a center of a maximum star.
A {\it pendant} tree is one for which each of its split vertices is adjacent to a leaf.
Here we show that all pendant trees $T$ have the leaf property for all $r\le \a(T)$.
We also consider pendant trees with exactly two split vertices and provide partial results on the locations of the centers of their maximum stars.
\end{abstract}

\medskip


\section{Introduction}

Let $[n]=\{1,\ldots,n\}$. Let $2^{[n]}$ and $\binom{[n]}{r}$ denote the family of all subsets and $r$-subsets of $[n]$ respectively. A family $\cF\subseteq 2^{[n]}$ is \textit{intersecting} if $F\cap G\neq \emptyset$ for $F, G\in \cF$. For any $\cF\subseteq 2^{[n]}$ and $x\in [n]$, let $\cF_x$ be all sets in $\cF$ that contain $x$. A classical result of Erd\H{o}s, Ko and Rado \cite{EKR} states that if $\cF\subseteq \binom{[n]}{r}$ is intersecting for $r\leq n/2$, then $|\cF|\leq \binom{n-1}{r-1}$. Moreover, if $r<n/2$, equality holds if and only if $\cF=\binom{[n]}{r}_x$ for some $x\in [n]$. This was shown as part of a stronger result by Hilton and Milner \cite{HilMil} which characterized the structure of the ``second-best'' intersecting families.

There have been multiple proofs of the Erd\H{o}s--Ko--Rado theorem. The original proof by Erd\H{o}s, Ko and Rado devised the now-central \textit{shifting} technique and used it in conjunction with an induction argument to prove the theorem. Daykin \cite{Dayk} demonstrated that the theorem is implied by the Kruskal--Katona theorem. Katona \cite{Katona} provided possibly the simplest and most elegant proof, a double counting argument using the method of cyclic permutations. More recently, Frankl--F\"uredi \cite{FraFur} gave another short proof that relied on a result of Katona on shadows of intersecting families, while we \cite{HurKamEKR} provided an injective proof using the aforementioned shifting technique. There have also been algebraic proofs, one using Delsarte's linear programming bound (see \cite{GodMea} and \cite{GodRoy} for details), and another using the method of linearly independent polynomials due to F\"uredi et al. \cite{Fur}.

The Erd\H{o}s--Ko--Rado theorem is one of the fundamental theorems in extremal combinatorics, and has been generalized in many directions. A very fine survey of the avenues of research, pursued as extensions of the Erd\H{o}s--Ko--Rado theorem, in the 1960's, 70's and 80's, is presented by Deza and Frankl \cite{DezFra}. In this note, we focus on a relatively recent graph-theoretic extension of the theorem.

\subsection{Erd\H{o}s--Ko--Rado graphs}

For a graph $G$ and integer $r\le \a(G)$, where $\a(G)$ is the size of a maximum-sized independent set in $G$, we define $\cI^r(G)$ to be the family of all independent sets of $G$ having size $r$.
For any family $\cF$ of subsets of $V(G)$ we denote by $\cF_x$ those sets of $\cF$ that contain the vertex $x$.
We call $\cI^r_x$ the {\it star centered on $x$}, and call $x$ the {\it star center} (use use the notation $\cI_x^r(G)$ in place of $\cI^r(G)_x)$. Call a graph $G$ $r$-EKR if, for any $\cF\sse I^r(G)$, $|\cF|\leq \textrm{max}_{x\in V(G)} |\cI^r_x(G)|$. 

Earlier results by Berge \cite{Berge}, Deza and Frankl \cite{DezFra}, and Bollobas and Leader \cite{BolLead}, while not explicitly stated in graph-theoretic terms, hint in this direction. The formulation was initially motivated by a conjecture of Holroyd, who asked if powers of the cycle graph on $n$ vertices are $r$-EKR for every $r\geq 1$. Holroyd's conjecture was later proved by Talbot \cite{Talbot}. The formulation also has connections with a fundamental conjecture of Chv\'atal \cite{Chvatal} on intersecting subfamilies of hereditary (closed under subsets) set systems. 

Holroyd and Talbot \cite{HolTal} made the following interesting conjecture about the EKR property of graphs. Let $\mu(G)$ be the size of a smallest \textit{maximal} independent set in $G$. 

\begin{cnj}
\label{j:minimaxconj}
For a graph $G$, let $1\leq r\leq \mu(G)/2$. Then $G$ is $r$-EKR.
\end{cnj}

Conjecture \ref{j:minimaxconj} appears hard to prove in general, but has been proven in a much more general form by Borg \cite{Borg} for $\mu(G)$ sufficiently large in terms of $r$.  In addition, it has been verified for certain graph classes. In the paper that introduced this graph-theoretic formulation of the EKR problem, Holroyd, Spencer and Talbot \cite{HolSpeTal} proved the conjecture for a disjoint union of complete graphs, paths and cycles containing at least one isolated vertex. Borg and Holroyd \cite{BorHol} later proved the conjecture for a certain class of interval graphs containing an isolated vertex. In \cite{HurKamChord}, we extended this result and verified the conjecture for all chordal graphs containing an isolated vertex. 

One of the reasons why verifying the conjecture for graph classes without isolated vertices is harder is that the intermediate problem of finding a center of a largest star is difficult. (It is easy to see that in a graph containing an isolated vertex, such a center is at an isolated vertex.) In this note, we consider this problem for trees. 
To that end, for a graph $G$ we define a vertex $x$ to be a {\it max} $r$-{\it center} if $|\cI^r_x(G)|$ is maximum among all stars of $G$.

In \cite{HurKamChord}, we proved that for any tree $T$ and $r\leq 4$, some leaf is a max $r$-center. 
The authors of \cite{EstrPast} call such trees $r$-{\it HK} (and {\it HK} if $r$-HK for all $r$).
We also conjectured that this is true for every $r\geq 1$. However, Baber \cite{Babe}, Borg \cite{BorCounter}, and Feghali, Johnson, and Thomas \cite{FegJohnTho} have separately shown that this conjecture is not true. 
This makes it interesting to consider trees for which the conjecture is true.

Define a vertex $v$ of a tree to be {\it split} if $\deg(v)\ge 3$.
The authors of \cite{FegJohnTho} consider a special class of trees called spiders: trees having exactly one split vertex.
(One can think of these as obtained from the star graph $K_{1,n}$, for some $n\geq 1$, by multiple subdivisions of edges.) They prove that two families of spiders, namely the family of all spiders obtained by subdividing each edge of the star graph exactly once (i.e. every leaf has distance two from the split vertex), and also the family of all spiders containing one leaf vertex adjacent to the split vertex, satisfy Conjecture \ref{j:minimaxconj}. Note that in both of these subfamilies of spiders, it is easy to find a vertex that is a max $r$-center (for any $r\geq 1$). 
In this note, we focus on the problem of determining the max $r$-centers in all trees with at most two split vertices. 
In Section \ref{s:SpiderStars} we prove (Theorem \ref{t:BestLeaf}) that some leaf of a spider is a max $r$-center and, in the process, also give a complete ordering on the sizes of all leaf stars.
In Section \ref{s:GenStars} we discuss trees with two split vertices and prove two theorems (Theorems \ref{t:pendant} and \ref{t:2split}) about the location of their max $r$-centers.

We first introduce some notation to describe spiders.

\subsection{Spiders}

Given a sequence of positive integers $L=(l_1,\ldots,l_k)$ we define the {\it spider} $S=S(L)$ to be the tree defined as follows.
The {\it head} of $S$ is the split vertex $v_0$ and, for $1\le i\le k$, the {\it leg} $S_i$ is the path $(v_0,v_{i,1},\ldots,v_{i,l_i})$.
We say that $L$ is in {\it spider order} if the following conditions hold:
\begin{enumerate}
\item
if $l_i$ and $l_j$ are both odd and $l_i<l_j$ then $i<j$,
\item
if $l_i$ and $l_j$ are both even and $l_i<l_j$ then $i>j$, and
\item
if $l_i$ is odd and $l_j$ is even then $i<j$.
\end{enumerate}

To simplify the notation somewhat, we will write $\cI^r_i(G)$ (respectively $\cI_{i,j}^r(G)$) in place of the more cumbersome $\cI^r_{v_i}(G)$ (respectively $\cI_{v_{i,j}}^r(G)$).

\section{Spider Star Centers}
\label{s:SpiderStars}

\begin{thm}
\label{t:Leafs}
Let $S=S(L)$ be a spider with $L=(l_1,\ldots,l_k)$ and suppose that $r\le\a(S)$.
Then for each $1\le i\le k$ and $1\le j<l_i$ we have $|\cI^r_{i,j}(S)|\le |\cI^r_{i,l_i}(S)|$.
\end{thm}

\begin{proof}
We define an injection $f:\cI^r_{i,j}(S)\rightarrow \cI^r_{i,l_i}(S)$.

Let $A\in \cI^r_{i,j}(S)$ and consider the path $P=(v_{i,j},\ldots,v_{i,l_i})$.
For $0\le h\le (l_i-j)$ we define $B$ by placing $v_{i,l_i-h}\in B$ if and only if $v_{i,j+h}\in A$; $B$ is the {\it flip} of $A$ on $P$, denoted $\flip_P(A)$.
Let $W=A-V(P)$; then set $f(A)=B\cup W$.

Clearly, $f(A)$ is independent, contains $v_{i,l_i}$, and has size $r$.
Also, if $f(A^\pr)=f(A)$, then $A^\pr=A$.
\end{proof}

\begin{thm}
\label{t:Head}
Let $S=S(L)$ be a spider with $L=(l_1,\ldots,l_k)$ and suppose that $r\le\a(S)$.
Then for every $1\le i\le k$ we have $|\cI^r_0(S)|\le |\cI^r_{i,l_i}(S)|$.
\end{thm}

\begin{proof}
For fixed $i$ we define an injection $f:\cI^r_0(S)\rightarrow \cI^r_{i,l_i}(S)$.

First we define $f$ to be the identity on $\cI^r_0(S)\cap \cI^r_{i,l_i}(S)$.

Second, let $A\in \cI^r_0(S)\setminus \cI_{i,l_i}^r(S)$ and consider the leg $S_i=(v_0,v_{i,1},\ldots,v_{i,l_i})$.
Write $v_{i,0}=v_0$ and, for $0\le h\le (l_i)$ we define $B$ by placing $v_{i,l_i-h}\in B$ if and only if $v_{i,h}\in A$; $B$ is the {\it flip} of $A$ on $S_i$, denoted $\flip_{S_i}(A)$.
Let $W=A-S_i$; then set $f(A)=B\cup W$.

Clearly, $f(A)$ is independent, contains $v_{i,l_i}$, and has size $r$.
Also, if $f(A^\pr)=f(A)$, then $A^\pr=A$.
\end{proof}

Together, Theorems \ref{t:Leafs} and \ref{t:Head} verify that for the family of spiders, max $r$-centers occur at leaves. In what follows, we not only find the best leaf of a spider but give a complete ordering of its leaves according to star size.

\begin{thm}
\label{t:BestLeaf}
Let $S=S(L)$ be a spider with $L=(l_1,\ldots,l_k)$ in spider order and suppose that $r\le\a(S)$.
Then for each $1\le i<j\le k$ we have $|\cI^r_{i,l_i}(S)|\ge |\cI^r_{j,l_j}(S)|$.
\end{thm}

\begin{proof}
Fix $i$ and $j$ with $1\le i< j\le k$.
We define an injection $f:\cI^r_{j,l_j}(S)\rightarrow \cI^r_{i,l_i}(S)$.
There will be three cases to consider, depending on the parities of $l_i$ and $l_j$.
By symmetry, we may assume that $l_i\not= l_j$.
First, we develop some terminology.

For a set $A\in\cI^r(S)$ we can define its ladder as follows.
A pair of vertices $\{v_{i,h},v_{j,h}\}$ ($1\le h\le\min(l_i,l_j)$) is called a {\it rung}, which we say is {\it odd} or {\it even} according to the parity of $h$.
A rung is {\it full} if both its vertices are in $A$.
The {\it ladder} $\cL$ of $A$ is the set of either all even or all odd rungs, depending on whether $v_0\in A$ or not, respectively.
$\cL$ is {\it full} if all its rungs are full.
If $\cL$ is not full then there is a first (i.e. closest to $v_0$) non-full rung $R$.
The partial ladder $\cL^\pr$ is the set of all (necessarily full) rungs above $R$.
Let $W$ denote the set of vertices of $A$ not on $S_i\cup S_j$.

First, we define $f$ to be the identity on $\cI^r_{i,l_i}(S)\cap \cI^r_{j,l_j}(S)$.

Next, we define the function $f$ on the remaining sets $A\in \cI^r_{j,l_j}(S)$ having partial ladders.
Define the path $P$ from $v_{j,l_j}$, up its leg to $R$, across $R$, and down the other leg to $v_{i,l_i}$; i.e. $P=(v_{j,l_j},\ldots,v_{j,h},v_{i,h},\ldots,v_{i,l_i})$, where $R=(v_{i,h},v_{j,h})$.
Now slide $A$ along $P$ until it contains $v_{i,l_i}$ --- the result we call $\slide_P(A)$.
Then set $f(A) = \cL^\pr \cup \slide_P(A) \cup W$.
Of course $|f(A)| = |A|$, $v_{i,l_i}\in f(A)$, $v_{j,l_j}\not\in f(A)$, and $f(A)$ is independent because $R$ was not full.
Moreover, $\cL^\pr(f(A)) = \cL^\pr(A)$, and so the inverse of $f$ on $f(A)$ is uniquely determined.

Note that in these first two cases $f$ preserves both inclusion and exclusion of $v_0$.
This means that $W$ cannot affect the independence of $f(A)$.

Finally, we define $f$ on the remaining sets $A$ having full ladders.
Spider order implies either that $l_j$ is even and $l_j<l_i$ or that $l_i$ is odd and $l_i<l_j$.
Having a full ladder implies that the former case has $v_{j,l_j}\in A$, while the latter case has $v_{i,l_i}\not\in A$.
In both cases these imply that $v_0\in A$.

When $l_j<l_i$ we let $P$ be the $v_{j,l_j}v_{i,l_j-1}$-path in $S$ (i.e. $P=v_{j,l_j},\ldots,v_{j,1},v_0,$ $v_{i,1},\ldots,v_{i,l_j-1}$).
When $l_j>l_i$ we let $P$ be the $v_{j,l_i-1}v_{i,l_i}$-path in $S$ (i.e. $P=(v_{j,l_i-1},\ldots,v_{j,1},v_0,$ $v_{i,1},\ldots,v_{i,l_i}$)).
In both cases we let $Q$ be the $v_{j,l_j}v_{i,l_i}$-path in $S$, minus $P$.
We shift $A$ along $P$ just one step toward $v_{i,l_i}$ --- call the result $\shift_P(A)$ --- and flip $A$ on $Q$ (that is, if $Q = (q_0, ..., q_k)$ then replace each $q_h$ in $A$ by $q_{k-h}$) --- call the result $\flip_Q(A)$.
Now define $f(A) = \shift_P(A) \cup \flip_Q(A) \cup W$.
Of course $|f(A)| = |A|$ and, because of the flip if $l_j<l_i$ or the shift if $l_j>l_i$, $v_{i,l_i}\in f(A)$, $v_{j,l_j}\not\in f(A)$, and A is independent.
Moreover, $f(A)$ has a full ladder, and so the inverse of $f$ on $f(A)$ is uniquely determined.

Notice that, because of the shift, $v_0\not\in f(A)$, and so $W$ cannot affect the independence of $f(A)$.
Thus the injection is complete.
\end{proof}

\section{General Star Centers}
\label{s:GenStars}

Call a spider {\it short} if one of its legs has length 1 ($l_1=1$ in spider order).
Given a set $D$ of vertices of a graph $H$, define a {\it (short) $D$-spidering} of $H$ to be any graph formed by adding a (short) spider to each vertex in $D$; that is, for each $v\in D$ add a (short) spider with center $v$.

A {\it thread} $P=v_1\cdots v_k$ in a graph $H$ is a path in which each interior vertex has degree two in $G$ ($\deg(v_i)=2$ for $1<i<k$).
An edge is one example of a thread.
For a set $D$ of vertices of a graph $H$, define $D$ to be a {\it thread-dominating set} if, for all $v\in V(H)−D$, $v$ has a thread to some vertex in $D$.

The authors of \cite{EstrPast} recently used the above flip and slide operations to prove the following theorem.

\begin{thm}
\label{t:ThreadDom}
If a graph $G$ is a short $D$-spidering of some graph $H$ with thread-dominating set $D$ then, for all $r\le \a(G)$, some leaf of $G$ is a max $r$-center.
\end{thm}

We define a tree to be {\it pendant} if every split vertex is adjacent to a leaf.
The following is then a corollary to Theorem \ref{t:ThreadDom}.

\begin{cor}
\label{c:pendant}
If $T$ is a pendant tree and $r\le\a(T)$ then some leaf off $T$ is a max $r$-center.
\end{cor}

\begin{proof}
The set of split vertices of $T$ is a thread-dominating set of $T$.
\end{proof}

For a leaf $x$ adjacent to vertex $y$ in a tree $T$, let $d'(x)=\deg(y)$.
For a graph $G$ that is a short $D$-spidering of some graph $H$ with thread-dominating set $D$, consider the following $r$-{\it Max Neighbor Property}: 
some max $r$-center is a leaf $x$ having maximum $d'(x)$ among all leaves.

Based on initial calculations, we imagined that it might be possible that, if a graph $G$ is a short $D$-spidering of some graph $H$ with thread-dominating set $D$ then $G$ has the $r$-max neighbor property for all $r\le\a(G)/2$.
Theorem \ref{t:BestLeaf} shows this to be true if $G$ is a short spider; i.e. a pendant tree with exactly one split vertex. 
We explored this possibility for the next simplest case: pendant trees with exactly two split vertices.
Below we present two theorems, the first showing that the $r$-max neighbor property fails in general, and the second showing a class of trees for which the $r$-max neighbor property holds.

We begin by defining a class of counterexample trees $T_k$ on $3k+7$ vertices, having $\a(T_k)=2k+4$ and $\mu(T_k)=k+3$.
Let $T_k$ have the path $(y,a,b,c,z)$ (which we call the {\it spine}), with $a$ having additional neighbors $u_1,\ldots,u_{k+1}$, $b$ having additional neighbors $v_1,\ldots v_k$, and each $u_i$ having the additional neighbor $w_i$ (see Figure \ref{f:fail}).
Then $T_k$ is pendant, and so Corollary \ref{c:pendant} says that, for all $r\le k+2$, some leaf is a max $r$-center.
The $r$-max neighbor property would say that $y$ is that leaf in $T_k$, which is true for $r\le 4$ but fails for $r\ge 5$.

\begin{figure}
\begin{center}
\begin{tikzpicture}[scale=.7]
\tikzstyle{every node}=[draw,circle,fill=black,minimum size=1pt,inner sep=2pt]
\draw node (y) [label=above: $y$] at (-4,0) {};
\draw node (a) [label=above: $a$] at (-2,0) {};
\draw node (b) [label=above: $b$] at (0,0) {};
\draw node (c) [label=above: $c$] at (2,0) {};
\draw node (z) [label=above: $z$] at (4,0) {};
\draw (y) -- (a) -- (b) -- (c) -- (z);
\draw node (u1) [label=left: $u_1$] at (-3.5,-2) {};
\draw node (w1) [label=left: $w_1$] at (-3.5,-4) {};
\draw node (u2) [label=left: $u_2$] at (-2,-2) {};
\draw node (w2) [label=left: $w_2$] at (-2,-4) {};
\draw node (u3) [label=left: $u_3$] at (-.5,-2) {};
\draw node (w3) [label=left: $w_3$] at (-.5,-4) {};
\draw (a) -- (u1) -- (w1);
\draw (a) -- (u2) -- (w2);
\draw (a) -- (u3) -- (w3);
\draw node (v1) [label=right: $v_1$] at (1,-2) {};
\draw node (v2) [label=right: $v_2$] at (3,-2) {};
\draw (c) -- (v1);
\draw (c) -- (v2);
\end{tikzpicture}
\end{center}
\caption{The tree $T_k$ with $k=2$.} 
\label{f:fail}
\end{figure}

\begin{thm}
\label{t:pendant}
For the tree $T_k$ we have $|\cI_z^r(G)|>|\cI_y^r(G)|$ for all $5\le r<\a(T_k)$.
\end{thm}

\begin{proof}
Let $X=\cI_z^r(G)\cap\cI_y^r(G)$, $Y=\cI_y^r(G)-X$, and $Z=\cI_z^r(G)-X$.
The set of elements of $Y$ not containing $c$ is in one-to-one correspondence with the set of those elements of $Z$ not containing $a$ --- given by the flip along the spine $(y,a,b,c,z)$.
Define $Y_c$ to be the remaining elements of $Y$ (these all contain $c$), and $Z_a$ to be the remaining elements of $Z$ (these all contain $a$).

Furthermore, the set of elements of $Y_c$ containing no $u_i$ is in one-to-one correspondence with the set of elements of $Z_a$ containing no $v_i$ --- again given by the flip along the spine $(y,a,b,c,z)$.
Define $Y'_c$ to be the remaining elements of $Y_c$ (these all contain some $u_i$), and $Z'_a$ to be the remaining elements of $Z_a$ (these all contain some $v_i$).

Let $t=r-2\ge 3$.
It is straightforward to calculate that $|Z'_a| = \binom{2k+1}{t}-\binom{k+1}{t}$, while $|Y'_c| = \binom{k+1}{t}(2^t-1)$.
Indeed, we show that $\binom{2k+1}{t}>\binom{k+1}{t}2^t$ for $t\ge 3$.
After canceling and clearing the denominators, the left side equals $(2k+1)(2k)(2k-1)m=(4k^2-1)(2k)m$, while the right side equals $(2k+2)(2k)(2k-2)m'=(4k^2-4)(2k)m'$, where the descending terms (if $t>3$) of $m$ are each greater, in turn, than those of $m'$.
Hence $|\cI_z^r(G)|>|\cI_y^r(G)|$.
\end{proof}

Next we define a restricted family of pendant trees having exactly two split vertices, for which the $r$-max neighbor property holds.
Let $T$ be such a tree, with $v_1$ and $v_2$ the two split vertices of $T$, labeled so that $\deg(v_1)>\deg(v_2)$, and choose leaves $v'_i$ adjacent to $v_i$ for each $i$.
Define the {\it spine} $S$ of $T$ to be the unique $v'_1v'_2$-path, with $S'=S-\{v_1,v_2\}$.
Denote by $T_i$ the subtree of $T-S'$ containing $v_i$ and observe that each $T_i$ is a spider, with $T_1$ having more legs than $T_2$.

Let $(l_{1,1},\ldots,l_{1,j_1})$ be the lengths of the spider $T_1$, and $(l_{2,1},\ldots,l_{2,j_2})$ be the lengths of the spider $T_2$, each written in spider order; recall that $j_2<j_1$.
We say that $T$ has a {\it spider embedding} if, for each $1\le j\le j_2$, the pair $(l_{1,j},l_{2,j})$ is in spider order.

\begin{thm}
\label{t:2split}
If $T$ is a pendant tree with exactly two split vertices, $r\le\a(T)$, and $T$ has a spider embedding, then $T$ has the $r$-max neighbor property; i.e. the leaf $v'_1$ is a max $r$-center.
\end{thm}

\begin{proof}
We label $T$ as above.
We prove this theorem by providing an injection $f$ from $\cI_2^r(T)$ to $\cI_1^r(T)$, where $\cI_i$ is shorthand for $\cI_{v'_i}$.
We partition the definition of $f$ according to how it acts on $S$, $T_1$, and $T_2$ (of course, each $v_i$ is in $S$ and $T_i$, but these separate definitions will agree on them).
Thus, for $X\in\cI(T)$ we write $X_S=X\cap S$ and $X_i=X\cap T_i$.

If $A\in \cI_1^r(T)\cap \cI_2^r(T)$ then $f(A)=A$.
In all other cases we have that $f$ flips $S$.
That is, for $A\in\cI_2^r(T)$, we have $f(A_S)=\flip(A_S)\ni v'_1$. 
If $\flip(A)\in\cI_2^r(T)$ for $A\in\cI_1^r(T)$ then $f=\flip$ on such $A$.
Otherwise, we further isolate the definition of $f$ along the legs $L_{i,1}$, $\ldots$, $L_{i,j_i}$ of $T_i$, with --- outside of exceptions that will be explained later --- $f(A_{i,j})\subset L_{3-i,j}$ for each $i$ and each $1\le j\le j_2$, where $A_{i,j}=A\cap L_{i,j}$.

Label the vertices of $L_{i,j}$ in adjacent order $(v_{i,j,0}$, $\ldots$, $v_{i,j,l_{i,j}})$, where $v_{i,j,l_0}=v_i$ and $v_{i,j,l_{i,j}}$ is the leaf.
Define $h_i$ to be the smallest integer such that $v_{i,j,h_i}\not\in A$ for each $i$.
Then, provided $h_i$ exists, $f$ will swap the elements of $A$ on $L_{1,j}$ and $L_{2,j}$ that ``precede'' $h_i$; that is, for each $k<h_i$, we place $v_{2-i,j,k}\in f(A_{i,j})$ if and only if $v_{i,j,k}\in A_{i,j}$.

At this point, all definitions are invertible whenever $f(A)$ is independent.
Thus we continue by extending $f$, modifying it to fix the cases in which $f(A)$ as currently defined is not independent.
Thus we may assume that $h$ does not exist.
In this case we swap the elements of $A$ on $L_{1,j}$ and $L_{2,j}$ up to $v_{i,j,k^*}$, where $k^*=\min\{l_{1,j},l_{2,j}\}$.
Observe that, since $l_{1,j}$ precedes $l_{2,j}$ in spider order, $f(A)$ is independent.
Indeed, if $l_{1,j}=l_{2,j}$ this is clearly true.
If $k^*=l_{1,j}<l_{2,j}$ then spider order implies that $l_{1,j}$ is odd, which means that $v_{1,j,l_{1,j}}\notin A$; hence $v_{2,j,l_{1,j}}\notin f(A)$ and so $f(A)$ is independent.
Finally, if $k^*=l_{2,j}<l_{1,j}$ then spider order implies that $l_{2,j}$ is even, which means that $v_{2,j,l_{2,j}}\notin A$; hence $v_{1,j,l_{2,j}}\notin f(A)$ and so $f(A)$ is independent.

This verifies the injection and therefore concludes the proof.
\end{proof}

\section{Open questions}

Determining whether or not spider graphs satisfy Conjecture \ref{j:minimaxconj} remains open. 
The compression/induction technique that has been used to prove Conjecture \ref{j:minimaxconj} for other graph classes appears difficult to use in this case. 
The nature of Theorem \ref{t:BestLeaf} implies that the max $r$-center may ``jump'' when we consider subtrees of the spider. 

In general, determining the max $r$-centers in trees, characterizing when such centers are leaves, and identifying which leaves they are (even for pendant trees with exactly two split vertices), all remain open problems.
In particular, the authors of \cite{EstrPast} pose that split vertices of trees are never max $r$-centers.
Thus we propose the following conjecture.
\begin{cnj}
If $T$ is a tree with no vertex of degree two then $T$ is HK.
\end{cnj}

Among several other interesting problems, they also propose finding which graphs have minimum-degree vertices as max $r$-centers, and ask whether or not every graph $G$ has a vertex of degree $\delta(G)+o(|G|)$ as the max $r$-center.

\end{document}